# Laboratorio di matematica:
# una sintesi di contenuti e metodologie

*Maria Dedò*
*Simonetta Di Sieno*

**Sunto**

In questo articolo ci proponiamo di tirare le fila di dieci anni di sperimentazione, da parte del centro "matematita", della modalità laboratoriale a tutti i livelli dell'insegnamento / apprendimento della matematica, analizzandone le potenzialità dal punto di vista degli studenti e dei docenti.

**Parole chiave**

Matematica, laboratorio, apprendimento attivo

**Introduzione**

La diffusa valutazione di inefficacia dell'insegnamento della matematica nella scuola italiana non è certamente andata affievolendosi negli ultimi anni. Non si tratta soltanto dei risultati di indagini più o meno accurate che confrontano gli esiti di tale insegnamento con quelli ottenuti in altri paesi, ma anche della sensazione di impotenza che assale molti insegnanti quando provano a conciliare aspettative proprie e aspettative degli studenti, richieste dal mondo del lavoro e *status* disciplinare/immagine della matematica, tempi rilassati adeguati a un apprendimento duraturo e ritmo imposto dall'istituzione.

La società in cui gli studenti si troveranno ad operare usciti da scuola sarà innegabilmente una società complessa con una grande varietà di processi che si susseguiranno in maniera estremamente rapida e quindi, per non essere solo vittime passive e indifese del cambiamento, i ragazzi a scuola hanno bisogno di acquisire strumenti di comprensione, di interpretazione della realtà e di conseguente intervento nella vita quotidiana. Così, non solo sono fortemente mutate le aspettative della società nei confronti delle competenze matematiche di chi esce dalla scuola, ma queste sono per di più difficilmente individuabili a priori, perché professioni diverse avranno bisogno di competenze matematiche di tipo diverso e, in alcuni casi, non si può nemmeno prevedere ora come queste muteranno nel giro di pochi anni. È chiaro allora che l'insegnamento della matematica dovrebbe adeguarsi a queste mutate condizioni, e che tale adeguamento richiede non tanto di sostituire un capitolo di matematica ad un altro, ma piuttosto, come dice Kra (in [K]), di "… esporre gli studenti al più vasto possibile panorama di fatti matematici" e, soprattutto, ai metodi del "far matematica".

E, per far ciò, il punto davvero cruciale, quello che "fa la differenza", sta negli insegnanti (non è un caso che il titolo del citato articolo di Kra sia proprio "Teachers are the key"). Al docente di matematica spetta il compito di progettare e realizzare, in stretta collaborazione con i colleghi ma in una maniera che sfrutti al meglio la specificità della sua competenza, percorsi formativi/educativi che promuovano negli alunni facoltà di ragionamento, pensiero razionale e capacità critiche, sviluppino la capacità di affrontare e risolvere problemi, sollecitino dimensioni operative e stimolino l'autonomia, la creatività personale e l'assunzione di responsabilità.

In questo quadro, il carattere di approntamento di una strategia che è tipico della matematica rappresenta un importante strumento di lavoro e una chiave di accesso ad altri saperi. Per ottenere risultati che non svaniscano nel giro delle poche settimane che passano tra una verifica dell'apprendimento e l'altra, i docenti devono compiere scelte forti che riguardano il *corpus* della disciplina, ma anche scelte significative che riguardano l'ambito delle relazioni interpersonali: essi hanno infatti la necessità – per garantire efficacia al loro lavoro – di sintonizzarsi con le motivazioni, gli interessi e le conoscenze pregresse degli studenti e compiere così anche scelte coerenti sul fronte della comunicazione e della gestione del gruppo classe.

L'esperienza degli ultimi dieci anni di sperimentazione del Centro "matematita" sul tema dei laboratori mostra che, in questo contesto, introdurre la modalità laboratoriale nell'insegnamento pre-universitario rappresenta una possibilità interessante; e, in questo articolo, ci proponiamo per l'appunto di tirare le fila dell'attività di questo decennio.

La modalità laboratoriale prevede che i ragazzi facciano esperienze di soluzione di problemi – nati dalla vita quotidiana o dall'interno della matematica stessa – utilizzando le conoscenze, matematiche ma non solo, trasmesse dalla scuola, collegando conoscenze acquisite in ambiti diversi e mostrando quindi quali fra di esse siano diventate davvero loro patrimonio personale (il che è naturalmente ben diverso rispetto al semplice mostrare di saper ripetere una data affermazione). Un insegnamento della matematica imperniato su attività di tipo laboratoriale permette ai ragazzi di lavorare insieme ai compagni con un obiettivo comune da raggiungere, esprimendo alti livelli di efficienza e di correttezza, consentendo la valorizzazione di abilità di solito sottovalutate o addirittura non riconosciute, e garantendo un maggiore coinvolgimento nell'apprendimento stesso, e di conseguenza una maggior efficacia di tale apprendimento, soprattutto sul lungo periodo.

Allo stesso tempo, i laboratori possono costituire per gli insegnanti una palestra di formazione e autoaggiornamento continuo, sia nella preparazione del laboratorio stesso, per la necessità di coniugare temi pregnanti e metodologie coinvolgenti, sia nella sua conduzione con i ragazzi e nella messa in atto e osservazione delle dinamiche costruttive cui si accennava e che, viceversa, ben difficilmente si ha modo di rilevare nell'arco di una lezione frontale.

Nel seguito di questo articolo, dopo aver dichiarato e analizzato cosa intendiamo per laboratorio, entreremo nei dettagli dei vari aspetti che abbiamo qui solo accennato, distinguendo la fase di preparazione (che coinvolge solo il docente), la fase di attuazione e la fase finale di valutazione del laboratorio; in particolare, relativamente alla seconda fase, discuteremo il ruolo che vi assumono diversi aspetti tipici della disciplina (dal rigore al linguaggio, dall'errore alla manipolazione di oggetti concreti o di animazioni virtuali) e, in particolare, la differenza di ruolo che questi rivestono in un'attività laboratoriale rispetto a una lezione frontale: mostrando così come non sia un caso che certi meccanismi potenzialmente positivi si manifestino più facilmente proprio in tale modalità di lavoro.

## Che cosa si intende per "laboratorio"?

Sono molte le accezioni con cui il termine "laboratorio" compare sia nella letteratura relativa all'insegnamento della matematica, sia nei testi di riferimento dei docenti, sia nella prassi scolastica corrente. Dato che non c'è unanimità su che cosa vada inteso come laboratorio, è utile qui premettere qual è la nostra "definizione" di laboratorio, che sarà quindi l'accezione del termine a cui faremo poi riferimento nel seguito di questo articolo.

Per indicare un'attività didattica con il termine "laboratorio", noi richiediamo che:
- gli studenti vi abbiano un ruolo attivo: non abbiano solo la parte di ascoltatori, ma debbano concretamente operare, lavorando a piccoli gruppi e discutendo fra di loro, per costruire le proprie conoscenze;
- gli insegnanti vi abbiano il ruolo della guida esperta che osserva e ascolta, che risponde a eventuali domande, che sa indirizzare su una via proficua e distogliere da una via poco significativa, e, soprattutto, che aiuta i ragazzi – alla fine – a tirare le fila dell'attività che hanno svolto.

Ci sono poi altri elementi che possono essere utili per individuare un laboratorio (ad esempio, l'utilizzo di materiale manipolabile, l'uso della tecnologia,...); ma, a ben vedere, questi non sono tanto elementi "necessari" per chiamare una certa attività con il nome di "laboratorio" quanto piuttosto elementi che possono in maniera abbastanza naturale portare a uno "spiazzamento" rispetto alla lezione frontale (in cui l'insegnante spiega e gli studenti ascoltano) e possono facilitare così una dinamica di tipo laboratoriale, quale quella che si è cercato di definire qui sopra.

Non è certo quella di questi anni la prima volta che, nella riflessione dei matematici e degli insegnanti di matematica, compare la modalità laboratoriale come opzione privilegiata dell'apprendimento-insegnamento. Basta pensare non solo a nomi relativamente recenti come quelli di Emma Castelnuovo o Vittorio Checcucci, ma anche a nomi ben più lontani nel tempo come quello di Giovanni Vailati (1863-1904), ricordando la sua battaglia per una scuola in cui gli studenti non siano costretti a "imparare delle teorie prima di conoscere i fatti a cui esse si riferiscono" o a "sentir ripetere delle parole prima di essere in possesso degli elementi sensibili e concreti da cui per astrazione si può ottenere il loro significato" (vedi [V]).

## Prima del laboratorio: la formazione dell'insegnante/animatore

La proposta del laboratorio si configura come un metodo di lavoro efficace per l'apprendimento della matematica, a qualunque livello scolastico ci si ponga, un metodo che permette di superare almeno in parte le presenti difficoltà e che offre una cifra unitaria all'esperienza matematica degli studenti. Tutto ciò richiede innanzi tutto una adeguata formazione degli insegnanti che si apprestano a proporlo ai loro studenti. Come scrive Giuliano Spirito in [Sp1], una didattica laboratoriale "prevede di necessità da parte dell'insegnante un di più di attenzione didattica, un di più di capacità di conduzione, persino un di più di chiarezza espositiva (quella chiarezza che discende da una meta-riflessione profonda sui nodi disciplinari). Ma in realtà una buona didattica laboratoriale richiede ancora e molto altro: richiede un ripensamento approfondito della gerarchia di contenuti che si vogliono veicolare e una riflessione puntuale – direi una riflessione problematica per problematica – su modi, tempi, nodi concettuali, ostacoli cognitivi dell'apprendimento. D'altra parte, solo per questa via diventa possibile far coesistere ambizioni alte (cioè, far conoscere agli alunni le acquisizioni fondanti del sapere matematico, quelle che rendono la matematica bella, prima ancora che utile) e efficacia didattica (cioè, ottenere che ciascun allievo raggiunga il massimo dei risultati a cui potenzialmente può accedere)."

Occorre quindi che i docenti vengano formati a tal fine, ponendo attenzione sia alla conduzione del gruppo classe che alla scelta dei concetti matematici coinvolti. Per le questioni del primo tipo (che riguardano sia il *setting* che il progetto pedagogico e che vanno dall'organizzazione e dall'uso dello spazio fino alla determinazione e all'uso del tempo così come dalla formulazione e negoziazione delle regole di comportamento fino alla individuazione dei ruoli), rimandiamo alla ricca bibliografia in proposito. In quel che segue invece ci preoccupiamo di descrivere come costruire (imparare a costruire, insegnare a costruire) percorsi matematici significativi da svolgere in laboratorio, ponendo attenzione alla pregnanza dei contenuti, alla individuazione di strumenti e sussidi per le attività di laboratorio e al rispetto di un consapevole filo conduttore.

Nell'esperienza del Centro "matematita", non sempre gli insegnanti/animatori hanno vissuto nel loro percorso scolastico (anche universitario, o post-universitario) esperienze di didattica per problemi (una metodologia di apprendimento che non appartiene solo alla matematica, ma che anzi è nata per i corsi di medicina: vedi [Sa]) o comunque non sempre sono sfuggiti al totalizzante metodo della lezione frontale: ciò rende necessaria un'attività di addestramento articolata che non è semplice né semplificabile. Così, innanzitutto, i partecipanti alla formazione sono invitati ad affrontare senza un'immediata precedente preparazione alcuni problemi e ad affrontarli insieme ad altre persone, il più delle volte in piccolo gruppo. In generale, i problemi devono essere di soluzione non immediata e non devono prevedere la necessità di ricorrere a tecniche sofisticate, ma piuttosto una rilettura (magari anche non banale) di risultati che possono essere anche molto "elementari"; devono soprattutto richiedere di usare quegli ingredienti di buon senso e corretto ragionamento che stanno alla base di qualsiasi acquisizione di lungo periodo di un concetto matematico (vedi [D3]).

Il fatto che i partecipanti abbiano un tempo sufficiente da dedicare a questi problemi e che possano confrontarsi nel piccolo gruppo porta quasi sempre alla risoluzione del problema affrontato. Ciò permette di discutere quali sono le condizioni ambientali, le relazioni fra allievi e docenti e quelle degli allievi fra di loro che favoriscono un'attività matematica comune e significativa (dotata di senso). Non solo. Questa attività di ricerca delle soluzioni permette di riprodurre in modo naturale le fasi in cui si articola la stessa attività del ricercatore matematico: il dover risolvere un problema di cui non si sa neppure se una soluzione sia possibile; il dover analizzare il problema "a buon senso" per inventarsi una soluzione, una via da seguire; il perdere tempo in ragionamenti inconcludenti e sbagliati, prima di trovare la strada giusta; il rileggere i risultati ottenuti anche ai fini di una loro generalizzazione (ovvero della produzione di nuovi "problemi"); il raccontare alla comunità scientifica i risultati trovati (vedi [C5]).

E ciò accade non solo nei gruppi di insegnanti! È davvero impressionante, nei resoconti che a volte i docenti ci mandano delle attività tenute in classe, constatare come quelli che l'insegnante segnala come i passaggi più significativi del lavoro dei suoi studenti (anche bambini delle prime classi della scuola primaria) siano, *mutatis mutandis*, quelli che abbiamo qui sopra indicato, e cioè anche i passaggi più significativi del lavoro di un matematico impegnato in un problema di ricerca.

I problemi più adatti a questa fase di primo impatto dei docenti con la modalità laboratoriale sono quelli, per un motivo o per l'altro, più spiazzanti. Ad esempio, a volte proponiamo un problema dall'aspetto del tutto innocuo come il seguente: "Ci sono tre porte chiuse. Dietro ad una c'è una Ferrari, dietro ad ognuna delle altre due c'è una scatola di cioccolatini. Siete invitati a scegliere una porta e a portarvi a casa quello che è nascosto dietro di essa (e supponiamo che voi siate interessati all'auto più che ai dolci). Voi scegliete una porta e, prima che la apriate, un amico apre una delle altre due e vi fa vedere che dietro c'è una scatola di cioccolatini. Voi potete ancora scegliere quale porta indicare per portarvi a casa la Ferrari. Tenete quella che avevate scelto oppure la cambiate? Perché?"

Molto spesso i docenti ritengono del tutto naturale rispondere: "è indifferente cambiare o tenere", ma questa risposta è sbagliata. Perché è sbagliata? Di che cosa non si è tenuto conto? Spesso, chi prova ad analizzare la situazione, senza farsi distrarre dal contesto, si rende conto che in due casi su tre è meglio cambiare la scelta iniziale, ma ne rimane così stupito da far fatica a convincersene. E ciò porta in maniera naturale a riflettere su quanto sia importante l'effetto di spiazzamento di una soluzione che non "proviene dall'alto", ma che esige di essere dimostrata al di là di ogni dubbio. Se non ci si cimenta direttamente nella soluzione di qualche problema, si perde la capacità di valutare la complessità della questione trattata, e la difficoltà quindi del compito: tutto si appiattisce e tutte le soluzioni sono uguali. E ricordiamoci che una dimostrazione (a qualunque livello, dalle risposte ai "perché" dei bambini della scuola primaria ai primi tentativi di giustificazioni formali dei ragazzi più grandi) non serve proprio a nulla, se non si è prima creata l'esigenza che ci sia qualcosa da dimostrare! Quindi, la ricerca del problema "strano" e dell'effetto sorpresa NON ha certo lo scopo di mettere la matematica sul piano di un gioco di prestigio; ma, viceversa, vuole far leva proprio sulla sorpresa per spingere a chiedersi perché succede quel che succede e per creare quindi l'esigenza della dimostrazione.

Fare esperienza diretta dei meccanismi con i quali ci si confronta con un problema in una situazione effettivamente laboratoriale, con le difficoltà di lettura di una situazione, di costruzione delle ipotesi di risoluzione, e di comunicazione dei risultati ottenuti diventa così un buon modo per capire quale valore possa avere un'analoga esperienza per gli studenti. E l'attenzione si concentra in maniera naturale sul tentativo di costruire buone, analoghe occasioni di stupore, di spiazzamento, di sorpresa. Pescando in quali ambiti? Con quali rapporti con la programmazione curricolare? Toccando i temi del consueto programma scolastico, o temi "alternativi" che di norma non vengono esplorati? Immaginando percorsi per tutti gli studenti? O solo per quelli "bravi"? O solo per quelli in difficoltà? Con quali rischi/vantaggi per la costruzione di un linguaggio rigoroso da parte degli studenti? Sono tutte domande che è necessario che il docente si ponga per proporre un laboratorio coerente con le scelte generali che ha compiuto con la classe.

Quando gli aspiranti docenti/animatori hanno scelto un tema all'interno di un ambito disciplinare fissato, viene proposto loro di costruire attorno ad esso sessioni di laboratorio per una certa successione di classi. Risulta evidente fin dalle prime proposte concrete che le attività legate a temi molto circoscritti, di scarso spessore culturale, sono spesso o la riproposizione banale di problemi/esercizi già visti mille volte per i quali non è difficile immaginare che sarà inesistente o molto ridotto l'interesse dei ragazzi oppure sono, semplicemente, esercitazioni camuffate da domande. Anche il fatto che gli esercizi vengano presentati sotto la forma di giochi non modifica la situazione: per esempio, i giochi con le monete per imparare decimi e centesimi incuriosiscono i molto piccoli, ma subito vengono guardati con sufficienza dai più grandi… Invece, chi riesce a individuare un nodo centrale in matematica, una questione importante su cui lavorare ha meno difficoltà nel costruire attività di laboratorio sensate e coinvolgenti che spesso aprono fronti diversi con la conseguente ovvia necessità di decidere/scegliere se spalancarli davvero oppure no. La fondamentale unità della matematica interviene presto a stabilire gerarchie di importanza nelle attività che vogliono illustrare (ad esempio) i numeri razionali: cade qualunque speranza di separare i numeri dalle forme, la geometria dall'aritmetica… (vedi, ad esempio, [CD]).

Un'evidenza sperimentale che abbiamo ricavato dal lavoro sul laboratorio di questi dieci anni è proprio il fatto che la scelta di un tema per un laboratorio è questione che non può essere impostata pensando solo a un segmento scolastico, ma deve essere affrontata tenendo conto il più possibile del percorso scolastico globale degli studenti, dalle prime classi del primo ciclo al termine della scuola pre-universitaria.

E ciò non contrasta affatto (anzi!) con le caratteristiche della disciplina che ci interessa: in matematica esistono alcuni nodi concettuali con i quali si devono confrontare sia gli insegnanti della scuola primaria che quelli della secondaria e per i quali le forme e le metodologie di presentazione alle varie età si richiamano fortemente. Certamente diversi sono gli obiettivi che si possono porre a un preadolescente che affronta nella scuola media una certa questione matematica da quelli che si possono fissare per un ragazzo delle prime classi di scuola primaria. Ma unica è la necessità di condurre l'uno e l'altro a fare direttamente esperienza dei risultati e dei metodi della matematica. E molto simili sono le maniere di accompagnare gli studenti a fare questa esperienza. Gli stessi nodi concettuali sono sottesi anche all'insegnamento della matematica nelle scuole superiori, soprattutto nei primi bienni, in maniera differente secondo le specificità degli indirizzi e delle tipologie, ma in maniera importante e significativa per tutti, con difficoltà e incomprensioni molto simili. In questa fase non sembra utile proporre una differenziazione fra bienni di scuole a matematica forte e bienni di scuole a matematica debole: l'esperienza di questi anni con le scuole superiori ha mostrato come lo scollamento fra docenti e studenti di queste età sia del tutto analogo nei due tipi di scuole e spesso provenga da un'analoga accelerazione sulla strada della formalizzazione delle scienze matematiche a discapito di una comprensione accettabile dei loro risultati e dei loro metodi. Evidentemente le acquisizioni disciplinari sono diverse nelle scuole a matematica forte e in quelle a matematica debole, ma ciò avviene a partire da contenuti significativi dal punto di vista disciplinare e non da tecnicismi superati e/o del tutto inutili, ma tradizionalmente accreditati di valenze formative generali, peraltro inesistenti.

E l'attenzione ai contenuti significativi fa anche sì che possano a volte comparire, magari sullo sfondo e non come principali attori delle sessioni di laboratorio, anche dei problemi davvero "difficili": ma è interessante che i ragazzi abbiano occasione di scorgere che esistono anche problemi di questo livello, problemi che loro non sanno risolvere, ma che magari anche i loro insegnanti non sanno risolvere e addirittura che nessuno al mondo sa ancora risolvere. Si dà un'idea ben piatta (e falsa!) della disciplina se li si convince, implicitamente, che la matematica serva soltanto a risolvere problemi banali (come quelli degli esercizi tutti uguali di alcuni libri di testo)!

Tornando al percorso di formazione dei docenti/animatori, nell'ultima fase viene loro chiesto di affidare la proposta che hanno costruito alla sperimentazione di qualche collega e di sperimentarla loro stessi. Le indicazioni contenute nelle schede per i ragazzi, le ipotesi fatte sulle reazioni e sui comportamenti degli studenti, le reazioni dei colleghi, le difficoltà che questi ultimi incontrano nel lavoro sono, tutti, elementi importanti per rendersi conto di che cosa sia o debba essere una sessione di laboratorio ben riuscita. Non importa se poi i docenti si limiteranno a usare *kit* di laboratori preparati da altri, perché l'aver capito "come funziona" impedisce – l'esperienza ce lo insegna – un uso acritico delle varie proposte e accresce il desiderio di intervenire direttamente nella costruzione di integrazioni e/o di totali novità. E in questa seconda fase emergono con chiarezza tutta una serie di questioni, strettamente connesse fra di loro, che provengono dalla osservazione del lavoro dei ragazzi in laboratorio e che illustriamo qui di seguito.

## Durante il laboratorio: il ruolo del rigore

Il laboratorio appare proprio come l'apoteosi del rigore sostanziale e l'annullamento del rigore formale. Elio Fabri, parlando dell'insegnamento della Fisica, scrive in [F]: "… il criterio del rigore non è quello delle esatte definizioni nel primo capitolo di un libro. Rigore significa chiarezza nel significato dei singoli passi, significa dire esplicitamente che i concetti si precisano man mano che si procede, che la validità di principi e teorie si rafforza quando se ne vede tutta la portata, che non ci sono singole leggi dimostrate da singoli esperimenti, ma che tutta la costruzione si regge nel suo insieme e nel suo insieme trova conferma nei fatti. Ottenere che lo studente capisca e ricordi tutto questo è più importante delle singole nozioni, regole, dati sperimentali. Ciò porta del tempo, ma è tempo ben speso, anche se si deve sacrificare qualche parte delle trattazioni tradizionali. Per convincersene, basta avere l'onestà di chiedersi quanto di quello che si fa in un corso con pretese di completezza viene effettivamente ricordato, anche solo dopo un anno, dallo studente medio: si arriverà necessariamente alla conclusione che la completezza senza chiarezza di comprensione è fatica sprecata. Questo non vuol dire naturalmente che nozioni, regole, dati sperimentali non debbano essere conosciuti e impiegati: ma che debbono esserlo in vista di uno scopo ben preciso e non fine a se stessi. ..."

Anche dopo tanti anni dalla pubblicazione, le parole di Fabri sono ancora attuali non sappiamo se per la Fisica, ma di sicuro per la Matematica. Gli esiti spesso sconcertanti che si ottengono quando si propongono riflessioni agli studenti (anche studenti universitari, di facoltà scientifiche) sui temi fondamentali con cui dovrebbero avere maggior confidenza sono solo un segnale dello scollamento che esiste fra insegnamento e apprendimento, un segnale d'allarme per l'evidente mancanza di consapevolezza del lavoro compiuto da parte degli studenti.

È proprio il fatto di centrare l'attenzione sul problema da risolvere – specie se tale problema è coinvolgente, o per la sua stessa natura, o per il fatto di discuterlo e affrontarlo insieme ai compagni – che porta in modo naturale a concentrarsi sulle difficoltà reali, sul nodo concettuale intorno al quale il problema è costruito. Servirà poi anche la tecnica - ovviamente nel laboratorio essa rappresenta un momento necessario - e i ragazzi sono i primi ad accorgersi che, sapientemente usata, a volte è proprio ciò che permette di compiere un passo avanti nella comprensione e nella gestione di un problema. Non si tratta però mai di tecnica artificiosa, o fine a se stessa, di esercitazioni sempre uguali che sembrano pensate per "ammaestrare" un soggetto da circo più che per insegnare un concetto.

E ciò che dà forza a questa situazione, e che pone le basi perché l'apprendimento possa risultare duraturo, è il fatto che il rigore non viene mai separato dal significato: è ovviamente ben diverso per lo studente appropriarsi di una tecnica che in prima persona ha scelto come la più adatta per venire a capo di una certa situazione (con l'occasione rinfrescando il suo funzionamento, che magari non ricordava perfettamente), piuttosto che vedere di questa tecnica solo un utilizzo artificioso in esercizi o problemi *ad hoc.*

È proprio nell'attività laboratoriale che si possono sciogliere in maniera naturale alcune incongruenze che potrebbero altrimenti sconcertare i nostri studenti. Il rigore non è mai un concetto assoluto, ma è qualcosa a cui ci si avvicina per successive approssimazioni, il che significa fra l'altro che la risposta a una domanda del tipo "è o non è abbastanza rigoroso?" comporta un raffronto fra l'affermazione in questione e il contesto in cui si sta lavorando. Gioco forza noi comunichiamo ai ragazzi, consapevolmente o inconsapevolmente, delle informazioni circa il rigore diverse, e a volte addirittura contraddittorie: c'è il momento in cui occorrono tutte le virgole al loro posto, e c'è il momento in cui occorre capire, grosso modo, "come funziona". Non c'è niente di male in questa incoerenza, che effettivamente rispecchia i diversi stadi che esistono quando si cerca di capire, di acquisire, di sistemare, di comunicare la matematica. Occorre però che i "patti" siano chiari con i ragazzi, e che l'incoerenza sia esplicita. Se resta sotto traccia e non chiara, comunica solo sconcerto (quando addirittura non viene interpretata come qualcosa che penalizza gli errori dei ragazzi e assolve quelli dell'insegnante!).

Il laboratorio appare anche come un luogo privilegiato dove esplorare le possibilità e le difficoltà che sono legate al linguaggio. Sappiamo benissimo quanto il linguaggio tecnico della matematica possa costituire una difficoltà, a volte davvero insormontabile, per molti ragazzi; esistono qui delle difficoltà reali, che l'insegnante non può cancellare e ignorare, nel tentativo di facilitare la vita ai suoi ragazzi. Può però – questo sì! – (anzi, deve) cercare di arrivare al nocciolo di queste difficoltà, in modo che la vera fatica sia spesa sulle difficoltà vere e non su tecnicismi artificiosi. Ed è proprio per questi motivi che il laboratorio è un contesto particolarmente adatto ad affrontare tali problemi: l'attività laboratoriale porta in maniera naturale a concentrare l'attenzione sulle difficoltà sostanziali, di senso, e rende spontaneo il non farsi distrarre dalle difficoltà formali.

Spesso si usano termini di cui si dà per scontato il significato, ma altrettanto spesso ci si rende conto di quanto sia necessario ripetere i concetti più volte, soffermandosi sul significato e l'uso corrente del termine per confrontarlo in un secondo momento con il significato che il termine ha nel linguaggio specifico della matematica.

Spesso accade che in laboratorio si parta da un linguaggio familiare, con l'uso di termini propri del linguaggio comune, che sembrano poter ben descrivere anche la situazione matematica in cui ci si trova, ma poi si passi a un linguaggio rigoroso quando ci si rende conto che, per capirsi, per comunicare, è necessario stabilire il significato del termine in quel contesto.

Come scrivono in [LT2] Domenico Luminati e Italo Tamanini: "Nel caso della matematica, un ruolo fondamentale per la comprensione viene giocato dalla formalizzazione del linguaggio con cui la si esprime. L'intuito, l'immaginazione e il buon senso, senza un'adeguata dose di formalismo, possono facilmente condurre a prendere abbagli e a cadere in errore."

## Durante il laboratorio: il ruolo dell'errore

Il laboratorio è il regno della matematica informale, quella che si fa in tuta da lavoro, e non in giacca e cravatta; e anche il ruolo dell'errore è naturalmente molto diverso nel momento in cui stiamo lavorando insieme, per cercare di capire qualcosa, rispetto al momento in cui stiamo cercando di ripulire, sistemare, formalizzare, comunicare ciò che abbiamo capito. E non soltanto per il motivo ovvio che di errori se ne fanno di più nella prima fase che nella seconda, ma per un motivo più intrinseco, ossia che gli errori, nella fase di ricerca, servono; di più, sono preziosi!

E questo è un fatto che l'insegnante si trova subito ad affrontare quando studia la possibilità di far lavorare i ragazzi in laboratorio; più in generale, fin dall'inizio si scontra con il ruolo che si può utilmente assegnare al livello informale di apprendimento in matematica. Tutti sono d'accordo sul fatto che possa fornire un fondamentale prerequisito per qualunque successiva acquisizione più formalizzata di sapere, ma è solo questo? In realtà è necessario fare un passo oltre, e operare delle scelte successive che siano coerenti con queste premesse. Se riteniamo che il momento informale sia un primo passo utile e necessario, allora deve essere chiaro anche ai ragazzi, dal patto educativo, che l'errore (perlomeno in questa fase) costituisce un normale strumento di conoscenza, qualcosa che può capitare a un gruppo di studenti, ma anche all'insegnante che si mette in gioco con loro (o che fa finta di mettersi in gioco…), qualcosa che occorre mettere a fuoco per poterne sfruttare al massimo le potenzialità ai fini dell'apprendimento, un bene che ci è prezioso. Ed è fondamentale creare per chi impara uno spazio da cui sia bandita la paura di sbagliare, che è uno dei nemici più potenti dell'apprendimento (vedi, ad esempio, [DS] per un'analisi più articolata del tema dell'errore).

Ciò è per noi talmente fondamentale che potrebbe addirittura costituire una "definizione" di laboratorio: un tipo di attività dove sbagliare non solo è lecito, ma è necessario; dove l'errore non viene represso, non viene cancellato, ma viene incoraggiato e discusso, per poter approfittare al massimo dell'analisi che ci induce a fare sui meccanismi del ragionamento, nostro e dei nostri allievi. Il che ci porta anche a suggerire che l'attività di laboratorio NON debba essere sottoposta direttamente a valutazione, ma sia vista come base dell'apprendimento i cui esiti, se ciò è richiesto, verranno – questi sì – sottoposti a valutazione. Sembra questa l'unica maniera per non far cadere i ragazzi nel meccanismo perverso della paura di sbagliare che è il primo nemico di ogni apprendimento.

## Durante il laboratorio: il ruolo della discussione

Un'altra caratteristica dei laboratori è il fatto che i ragazzi, lavorando in gruppo, sono "costretti" a discutere gli uni con gli altri e a fare lo sforzo di estrinsecare ai compagni la loro visione di un dato problema. Si tratta di un punto fondamentale: come è stato ampiamente studiato, esiste un vero e proprio stacco tra il momento della comprensione di un dato concetto o di un dato problema e il momento in cui lo stesso concetto è da ritenersi acquisito, al punto da essere in grado di raccontarlo, a se stessi e agli altri.

Non si tratta solo di uno stacco temporale, quanto piuttosto di uno stacco sostanziale: tant'è vero che non è affatto automatico il passaggio dall'uno all'altro livello ed è quindi necessario che l'insegnante preveda delle attività *ad hoc* che aiutino ogni studente ad acquisire questa consapevolezza. La discussione all'interno di piccoli gruppi forza i ragazzi a cimentarsi in queste attività di consolidamento metacognitivo, che altrimenti potrebbero apparire inutili e anche un po' noiose a chi pensa: "tanto io ho già capito come funziona".

Va osservato che, al di fuori di una didattica laboratoriale, è davvero poco frequente che i ragazzi abbiano occasione di parlare di argomenti di matematica fra di loro: la comunicazione (su temi matematici) non è mai *inter pares*, ma è sempre comunicazione fra chi sa e chi non sa. E tutti gli aspetti che abbiamo ora esaminato (il rigore, il linguaggio, l'errore) hanno una valenza assai diversa nella comunicazione paritetica rispetto a quella che paritetica non è. In una comunicazione paritetica, occorre in primo luogo capirsi, ed è proprio questa necessità di capirsi a vicenda che porta a mettere l'accento sul rigore di sostanza e a cercare un linguaggio non ambiguo, che conduce a capire in maniera molto naturale che un errore (opportunamente discusso) può a volte farci fare dei passi avanti nella comprensione di un problema.

E, soprattutto, la discussione paritetica informale porta gli allievi a "capire di aver capito", se hanno capito. Non si tratta di un gioco di parole: quante volte i nostri studenti ci danno l'impressione di avere più o meno acquisito un dato concetto, ma… di non aver chiaro, loro stessi, il fatto che l'hanno acquisito! A volte (e questo aspetto è fortemente evidente in alcune situazioni, ad esempio con gli studenti del corso di laurea in Scienze della Formazione Primaria) il mestiere dell'insegnante consiste non tanto nel proporre concetti nuovi, quanto piuttosto nel rendere l'interlocutore consapevole di ciò che già sa.

E naturalmente, sotto questo aspetto, un momento fondamentale del laboratorio è quello finale in cui l'intero gruppo classe, con la guida dell'insegnante, tira le fila del lavoro svolto: è in questa fase che l'insegnante può compiere un'opera preziosa mettendo in evidenza quei momenti che hanno costituito dei passaggi significativi, discutendo le differenze fra gli approcci di diversi gruppi, in modo che i ragazzi acquisiscano in pratica l'idea che non c'è mai un'unica strada per risolvere un problema, mettendo in luce sia quel *quid* che ha fatto fare un salto qualitativo in direzione dell'astrazione, sia il ruolo della tecnica, ad esempio una buona notazione che abbia permesso di gestire in maniera tranquilla un problema che altrimenti sembrava estremamente difficile.

## Durante il laboratorio: concreto e astratto, reale e virtuale

Il laboratorio (anche facendo riferimento ai più comuni laboratori di fisica o di altre discipline scientifiche) è spesso identificato come un'altra stanza, diversa dalla normale aula scolastica, dove sono a disposizione oggetti da manipolare e dove si possono fare esperimenti. Dicevamo all'inizio come la presenza di oggetti di questo tipo possa essere utile, e come altresì non sia caratterizzante dell'attività laboratoriale per la matematica. L'uso di oggetti concreti al fine di rappresentare esempi dei problemi e dei concetti astratti della matematica è qualcosa di assai delicato, dalle enormi potenzialità, ma anche da usare con una certa cautela.

Intanto, non tutti i problemi e gli argomenti si prestano a una modellizzazione tramite oggetti concreti e, d'altra parte, questa risulta utile (e, anzi, preziosa) solo quando sorge in modo naturale da una contestualizzazione (realistica!) del problema: gli esempi artificiosi (per voler infilare il concreto a tutti i costi) risultano non solo inutili, ma anche dannosi. Abbiamo in mente alcuni libri di testo che propongono di trovare il volume di pentole a forma di prisma triangolare con lato di base di lunghezza 1 metro, o di portaombrelli che sono cilindri cavi, con cavità a forma di piramide quadrangolare...: e naturalmente, se questo deve essere il concreto, meglio restare sul piano astratto!

Inoltre, anche nei casi in cui una modellizzazione concreta si presta particolarmente bene a illustrare il problema che si sta discutendo, occorre non dimenticare che l'oggetto concreto non è mai il concetto astratto che esso vuole rappresentare: per quanto questa affermazione sia così banale, è opportuno a nostro avviso ribadirla, perché tendiamo (noi per primi!) a dimenticarcene; non perché non lo sappiamo, ma perché, più il concetto (esempio: la sfera) è chiaro nella nostra testa, più lo proiettiamo sull'oggetto concreto che lo rappresenta (esempio: un pallone da calcio), arrivando addirittura a non vederne le accezioni concrete (il materiale, il peso, il colore, le cuciture fra pentagoni ed esagoni che ne costituiscono le facce...), ma a vedere solo l'idealizzazione del concetto astratto che esso rappresenta. Il che ovviamente non può e non deve accadere a chi non conosce il concetto astratto che noi vorremmo con quell'oggetto concreto rappresentare (vedi [D2]).

È un po' come quando facciamo una foto perché una certa cosa (esempio: il volo di un gabbiano) ci ha colpito particolarmente e ci sembra bellissimo: non ci accorgiamo nemmeno che è troppo piccolo, o troppo lontano, e magari che in primo piano c'è, sgradevolissima, una discarica di immondizia; non ce ne accorgiamo perché nel nostro immaginario ciò che ci ha colpito era talmente bello e importante che riempiva completamente il campo visivo; ma quando confrontiamo il campo visivo mentale della nostra immaginazione con il campo visivo reale della foto scattata ci accorgiamo… che la foto è bruttissima, perché contiene tante cose che non avevamo… visto e che sono quelle che vede una persona "normale", che non aveva i nostri motivi per vedere soltanto il volo del gabbiano.

Va detto che fraintendimenti di questo genere (applicati alla comunicazione della matematica attraverso modelli concreti) accadono più facilmente nella comunicazione (NON paritetica) tra docente e studenti, ovvero nella comunicazione tra chi ha già in testa un'immagine mentale di riferimento del concetto astratto a cui il modello si riferisce e chi viceversa se la sta costruendo; è più difficile che analoghi fraintendimenti si verifichino nell'ambito della comunicazione paritetica che rappresenta la normalità in laboratorio.

Al docente resta principalmente da ricordarsi di questo problema nel momento in cui sceglie il materiale che accompagnerà il laboratorio. Ad esempio, noi sappiamo già che cos'è un segmento e quindi vediamo senza esitazione un segmento nell'oggetto fotografato qui sotto: capita addirittura che non ci accorgiamo proprio più delle caratteristiche per cui esso si guarda bene dall'essere un segmento. E dobbiamo invece ricordarcene quando decidiamo il materiale che può essere funzionale a un dato laboratorio: può essere sensato, ad esempio, ipotizzare di usare materiale di questo tipo con i ragazzi più grandi della scuola secondaria, ma privilegiare materiali diversi con i bambini della scuola primaria.

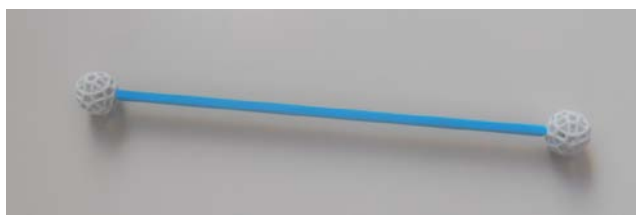

Infine, occorre tenere sempre ben presente che l'obiettivo dell'attività laboratoriale è l'acquisizione (più o meno solida a seconda dell'età e del tipo di problema) di un certo concetto astratto: occorre quindi stare attenti a che la manipolazione di oggetti concreti si configuri come uno strumento, volto a dare sostanza e significato a tale concetto astratto, e che non rischi invece di sostituirsi al concetto astratto stesso (vedi a questo proposito la polemica di Russo in [R]); d'altra parte, se si mantiene questa attenzione, è proprio l'acquisizione del concetto astratto che si giova della manipolazione di oggetti concreti che danno sostanza e significato all'astrazione che stiamo costruendo.

Pensiamo, per fare un esempio, alla dinamica tra astratto e concreto che si pone quando ci domandiamo (con ragazzi abbastanza grandi, della scuola secondaria) che cosa può significare la "lunghezza di un segmento", riferita a un materiale come quello raffigurato nella foto sulla destra. È una bella scoperta rendersi conto che il materiale (che è costruito in maniera molto raffinata, proprio dal punto di vista della matematica) presuppone che "lunghezza" di una data bacchetta significhi distanza fra i centri delle due sferette inserite agli estremi della bacchetta stessa; e che questa nozione di lunghezza è proprio quella che permette, disponendo di bacchette le cui lunghezze sono in determinati rapporti (p.es. rapporto aureo), di effettuare con gli oggetti concreti le costruzioni geometriche che si potrebbero fare con i corrispondenti segmenti astratti.

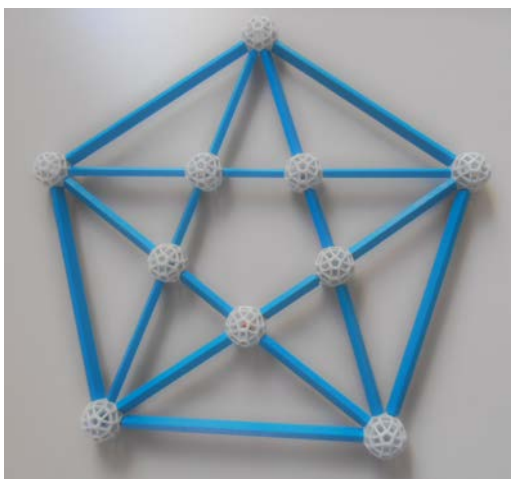 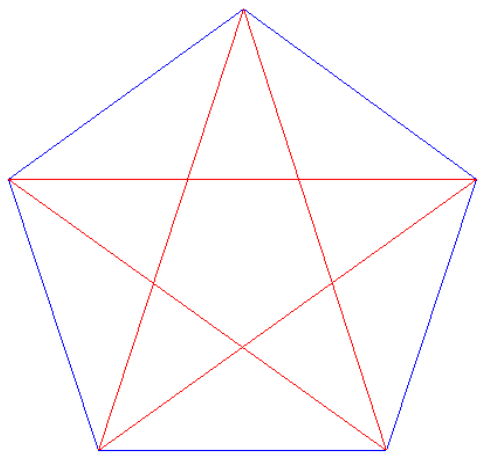

Anche l'interazione con ciò che la tecnologia ci mette oggi a disposizione è – come la manipolazione di oggetti concreti – un aspetto che può essere utile per un'attività laboratoriale, ma che non è caratterizzante. Ci possono essere laboratori che utilizzano la tecnologia e gli aspetti virtuali, e ci possono essere laboratori che non la utilizzano, così come esistono maniere di interagire in classe con la tecnologia che rappresentano attività di tipo laboratoriale e altre che di laboratoriale non hanno proprio nulla! Perché si possa parlare – l'abbiamo data come una specie di definizione – di attività laboratoriale, occorre che il ruolo dei ragazzi nell'interazione con il virtuale sia attivo e non passivo. Si potrebbe dire un "davvero" attivo: questo per sottolineare che non ci si accontenta del fatto che i ragazzi davanti a uno schermo possano fare un *clic*, dopodiché succede qualcosa di già predisposto che essi possono solo osservare. Piuttosto, è davvero interattivo un contesto in cui l'animazione predispone una sorta di ambiente di lavoro nel quale l'utente possa fare degli esperimenti in autonomia, e registrarne gli esiti. E meglio ancora, naturalmente, se i ragazzi possono essere loro in prima persona a realizzare degli oggetti virtuali.

## Dopo il laboratorio: la valutazione

Abbiamo già argomentato per quali motivi riteniamo che l'attività in laboratorio non debba essere oggetto di valutazione per il singolo ragazzo. Troppo spesso la scuola italiana è scuola che riduce il suo compito alla valutazione, mentre noi pensiamo che la scuola debba restare *in primis* una scuola che insegna!

Diverso invece è il discorso sulla valutazione del percorso laboratoriale compiuto con il gruppo classe; questa appare infatti come un elemento di grande importanza e che va affrontato, approntando strumenti che consentano di misurare l'efficacia di tale percorso. Senza entrare nel merito di come si possano costruire tali strumenti di valutazione, ci limitiamo qui a sottolineare alcune questioni da tenere presenti nella loro predisposizione.

Innanzitutto, una qualunque valutazione deve necessariamente prevedere una "misura" di efficacia sul medio-lungo periodo; appare infatti del tutto priva di senso la valutazione di un elemento troppo limitato, sia perché l'apprendimento su un segmento breve resta giocoforza un apprendimento fragile, sia perché tutto ciò che abbiamo messo in evidenza sulla situazione laboratoriale punta viceversa a una visione complessiva del percorso.

Un altro elemento che va messo in evidenza è l'opportunità (la necessità?) di un lavoro collettivo di discussione fra gli insegnanti per mettersi al corrente a vicenda di un'attività laboratoriale e per valutare l'efficacia del lavoro svolto: nella nostra esperienza accade infatti in maniera abbastanza frequente che una discussione di questo genere porti alcuni insegnanti a cambiare il proprio punto di vista, anche radicalmente, magari perché li porta a tener conto di alcuni aspetti più difficilmente quantificabili che avevano in prima battuta ignorato. Tanto per fare un esempio, si può citare il livello di responsabilizzazione in prima persona degli studenti, che è ovviamente un elemento cruciale perché l'apprendimento possa fissarsi e accrescersi nel tempo e che è sicuramente un aspetto privilegiato nelle attività laboratoriali.

## Conclusioni

Come dicevamo nell'introduzione, i temi discussi in questo articolo hanno origine da un'attività di dieci anni del centro "matematita" diretta alle scuole. Ci sembra utile a questo punto dare qualche numero di tale attività:

- i laboratori tenuti presso il Dipartimento di Matematica hanno toccato più di 1000 classi, quindi più di 20000 studenti e più di 2000 insegnanti
- i giochi *online* proposti dal sito http://www.quadernoaquadretti.it/ hanno toccato fra le 3000 e le 4000 classi di scuola primaria e secondaria di primo grado, quindi tra i 1500 e i 2000 docenti e circa 100000 ragazzi;
- i *kit* di laboratorio prestati alle scuole (l'attività ha avuto luogo solo negli ultimi 5 anni) hanno toccato più di 500 classi, quindi più di 10000 studenti e un migliaio di insegnanti.

Sono questi numeri (a cui si aggiunge l'esperienza di numerosi corsi di formazione per docenti tenuti in questi anni) che danno valore ai pareri positivi, e in alcuni casi addirittura entusiastici di tanti insegnanti che hanno sperimentato queste attività; e che ci hanno testimoniato, a distanza di anni, come i ragazzi che sono stati esposti a percorsi di questo tipo abbiano incontrato successivamente difficoltà assai minori rispetto ai loro compagni nelle prove di valutazione nazionali.

Ovviamente la modalità laboratoriale non è l'unica maniera per insegnare matematica, ma, in questo momento storico, essa appare come una maniera particolarmente efficace per farlo.

**Bibliografia**

[AAVV1] *I mosaici nell'insegnamento della geometria elementare*, a cura del gruppo per la scuola elementare Dip. di Matematica, Univ. Studi Milano,

[AAVV2] CD *matemilano, percorsi matematici in città*, Kangourou ed., Milano 2004, SIAE…

[B] M. Bertolini, G. Bini, P. Cereda, O. Locatelli, *Passeggiare tra le superfici*, collana "Quaderni di laboratorio", Mimesis 2012.

[Bo] G. Bolondi, *Metodologia e didattica: il laboratorio*, Rassegna, Periodico dell'Istituto pedagogico italiano, anno XIV, aprile 2006, pp. 59-63.

[BI] A. Brigaglia, G. Indovina, *Stelle, girandole e altri oggetti matematici*, Decibel ed., Padova 2003 ISBN 88-08-07857-4

[C1] M. Cazzola, *Per non perdere la bussola*, Decibel ed. Padova 2001 ISBN 88-08-07855-8

[C2] M. Cazzola, *Problem-based Learning and Mathematics: Possible Synergical Actions*, in L. Gòmez-Chova, D. Martì Belenguer, I. Candel Torres (editors), *Proceedings of ICERI 2008 Conference,* (Madrid, Novembre 2008), Valencia, IATED.

[C3] M. Cazzola, *Problem-based Learning and Teacher Training in Mathematics: the Role of the Problem*, in M. Tzekaki, M. Kaldrimidou, H. Sakonidis (editors), *Proceedings of the 33rd Conference of the Internetional Group for the Psychology of Mathematics Education* (Salonicco, Luglio 2009).

[C4] M. Cazzola, *WIMS all'Università di Milano-Bicocca,* Tecnologie didattiche, vol.19, n.3, pagg. 170-175.

[C5] M. Cazzola, *L'insegnamento della matematica: una didattica meta cognitiva,* in O. Albanese, P.A. Doudin, D. Martin (a cura di), *Metacognizione ed educazione,* Franco Angeli, Milano, 2003.

[CD] E. Colombo, M. Dedò, *Uguali? Diversi!,* collana- "Quaderni di laboratorio", Mimesis, in preparazione.

[D1] M. Dedò, *Rigour in Communicating Maths: a Mathematical Feature or an Unnecessary Pedantry?*, in E. Behrends, N. Crato, J.F. Rodrigues (editors), *Raising Public Awareness of Mathematics*, Springer, 2012, pagg. 339-358.

[D2] M. Dedò, *Matemática informal: ¿una contradicción?*, in R. Mallavibarrena (editor), *Escuela de Educación Matemática "Miguel de Guzmán": enseñar divulgando*, Secretaria General Tecnica, 2010. - ISBN 978-84-369-4937-7. - pp. 47-72

[D3] M. Dedò, *Più matematica per chi insegna matematica*, La Matematica nella Società e nella Cultura, *Bollettino dell'Unione Matematica Italiana,* (8), 4-A, Agosto 2001, pp. 247-275.

[DR] A. Deledicq, R. Raba, *Il mondo delle pavimentazioni*, ed. Mimesis – Kangourou Italia, Milano ISBN 88-8483-137-7

[DiS1] S. Di Sieno, *"Doing Mathematics": a Crucial Step in Mathematical Teachers Training,* in "Trends and Challenges in Mathematics Education" Jianpan Wang, Binyan Xu (eds), East China Normal University Press, pp. 365-375, Shangai,2004.

[DS] M. Dedò, L. Sferch, *Right or Wrong? That is the Question*, Notices of the Amer. Math. Soc. , vol. **59**, n.ro 7 (Agosto 2012), pagg. 924-932.
(http://www.ams.org/notices/201207/rtx120700924p.pdf).

[F] E. Fabri, *Che cos'è il rigore logico in fisica?,* in AAVV, *Guida al laboratorio di fisica,* Zanichelli, Bologna, 1995, http://www.df.unipi.it/~fabri/rigore/rigore00.htm

[K] I. Kra   *(Math) Teachers are the Key*, Notices AMS, vol. 59, n.4, aprile 2012

[Li] M. Liverani, *Qual è il problema?*, Mimesis ed., Milano 2005

[Lo] O. Locatelli, *Torri e serpenti*, collana "Quaderni di laboratorio", Mimesis 2006.

[LT1] D. Luminati, I. Tamanini, *Problemi di massimo e minimo*, collana "Quaderni di laboratorio", Mimesis 2009.

[LT2] D. Luminati, I. Tamanini, *La matematica a casa tua,* L'insegnamento della matematica e delle scienze integrate, vol. 28 A n.1, Gennaio 2005, pagg. 33-60.

[Ra] M. Rampichini, *Una palestra di matematica interattiva*, L'insegnamento della matematica e delle scienze integrate, vol. 35 A-B n.3, Maggio-Giugno 2012, pagg. 296-320.

[Ru] L. Russo, *Segmenti e bastoncini*, Feltrinelli, 1998.

[Sa] J.R.Savery, *Overview of Problem-based Learning: Definitions and Distinctions*, The Interdisciplinary Journal of Problem-based learning, vol.1, n.1, 2006.

[Sp1]   G. Spirito, *A proposito di metodologie, a proposito di contenuti*
in http://www.quadernoaquadretti.it/scuola/riflessioni/spirito_08.pdf

[Sp2] G. Spirito, *La lotteria di Babilonia*, Decibel ed., Padova 2003,  ISBN 88-08-07853-1

[V] G. Vailati, Recensione di *C.Laisant. La mathématique: philosophie, enseignement,* in "Giovanni Vailati, Scritti", a cura di M. Quaranta, Forni, Bologna, 1987.

[W]    H. Wu, *The Mis-education of Mathematics Teachers,* Notices of the Amer. Math. Soc., vol, 58, n. 3, March 2011. http://www.ams.org/notices/201103/index.html.